\documentclass[11pt]{amsart}
\usepackage{enumerate,amssymb,amsfonts,latexsym,psfrag}
\usepackage{amscd,amsmath}
\usepackage[dvips]{graphicx}

\newcommand{\cF}{{\mathcal{F}}}
\newcommand{\cS}{{\mathcal{S}}}

\newcommand{\cR}{{\mathcal{R}}}

\newcommand{\R}{\mathbb{R}}

\newcommand{\I}{^{-1}}

\newcommand{\hX}{\hat{X}}
\newcommand{\hf}{\hat{f}}

\newcommand{\sone}{{\mathbb{S}}^1}
\newcommand{\sph}{{\mathbb{S}}^2}
\newcommand{\diam}{\mbox{diam}}  
\newcommand{\ra}{\rightarrow}
\newcommand{\bs}{{\bar s}}

\renewcommand{\sb}{\overline{s}}

\newcommand{\gm}{$f\colon G\to G$\ }

\newcommand{\gmx}{$f\colon G\to G$}

\newcommand{\tgm}{$F\colon\ofg\to\ofg$\ }

\newcommand{\tgmx}{$F\colon\ofg\to\ofg$}

\newcommand{\trmx}{$\phi\colon\tau\to\tau$}

\newcommand{\trtrmx}{$\Phi\colon\cR\to\cR$}

\newcommand{\tPhi}{\tilde{\Phi}}

\newcommand{\ofg}{\mathbb{G}}
\newcommand{\ofi}{\mathbb{I}}
\newcommand{\oft}{\mathbb{T}}

\newcommand{\ofr}{\mathbb{R}}
\newcommand{\ofc}{\mathbb{C}}
\newcommand{\ofv}{\mathbb{V}}
\newcommand{\ofe}{\mathbb{E}}

\newcommand{\ol}[1]{\overline{#1}}

 \newcommand{\pichere}[2]
 {\begin{center}\includegraphics[width=#1\textwidth]{#2}\end{center}}
 \newcommand{\lab}[3]{\psfrag{#1}[#3]{$\scriptstyle{#2}$}}

\newcommand{\Int}{\mathrm{Int}}

\renewcommand{\mod}{\mathrm{mod}}

\newtheorem{thm}{Theorem}[section]

\newtheorem{lemma}[thm]{Lemma}
\newtheorem{prop}[thm]{Proposition}

\theoremstyle{definition}

\newtheorem{defn}{Definition}
\newtheorem{defns}[defn]{Definitions}
\newtheorem{example}{Example}
\newtheorem{examples}[example]{Examples}

\theoremstyle{remark}

\newtheorem{rmk}{Remark}
\newtheorem{rmks}[rmk]{Remarks}

\title{Extensions, quotients and generalized pseudo-Anosov maps}
\dedicatory{Dedicated to Dennis Sullivan on his
$60^{\text{th}}$ birthday}
\author{Andr\'e de Carvalho}
\address{Institute for Mathematical Sciences, 
State University of New York at Stony
Brook, NY 11794-3660, USA}
\email{andre@math.sunysb.edu}
\date{August 2001}

\begin{document}

\begin{abstract}
We describe a circle of ideas relating the dynamics of 2-dimensional
homeomorphisms to that of 1-dimensional endomorphisms. This is used to
introduce a new class of maps generalizing that of Thurston's
pseudo-Anosov homeomorphisms.
\end{abstract}

\maketitle

\bibliographystyle{amsalpha}

\section{Introduction}
In this paper we discuss a circle of ideas which is present in many
different contexts in dynamical systems. It was first
introduced by Williams in his study of expanding attractors, and
has been used since by many authors. In its most basic form
it can be stated as follows: 

\smallskip
\noindent
{\em Collapsing segments of stable manifolds of a homeomorphism yields
a lower dimensional endomorphism; the original homeomorphism may be
recovered by taking the natural extension (inverse limit) of the
quotient endomorphism.}

\smallskip
The context we focus attention on is the interplay this creates
between the dynamics of 1-dimensional endomorphisms --- endomorphisms
of trees and graphs --- and that of 2-dimensional homeomorphisms ---
homeomorphisms of surfaces. In recent years, this interplay between
graph endomorphisms and surface homeomorphisms has been used for
different, although related, purposes. For example, it was used to
give algorithmic proofs of Thurston's classification theorem for surface
homeomorphisms up to isotopy. It also appeared in the contruction of
models for families of surface homeomorphisms passing from trivial
to chaotic dynamics as parameters are varied. And a combination of these
results led to a conjecture about the way in which {\em forcing} organizes 
the braid types of horseshoe periodic orbits. 

Complexification should yield a closely related discussion --- which
will not be treated here --- linking the dynamics of endomorphisms of
branched (Riemann) surfaces to the dynamics of automorphisms of
$\ofc\,^2$. The natural extension of an endomorphism of a branched
surface is usually a homeomorphism of a surface lamination. Such
laminated spaces have already been found to be among the main objects
in the study of complex H\'enon maps. In both the real and complex
cases, much is known about 1-dimensional endomorphisms, at least when
the space is unbranched, whereas about homeomorphisms in dimension 2
much less is known. 

This double interplay --- between 1- and 2-dimensional dynamics and
between the real and complex settings --- seems to be an interesting
approach to explore. In this paper we develop some
of its aspects, mostly in the real setting.

We also introduce a different kind of quotient, taking a surface
homeomorphism to another surface homeomorphism by collapsing
dynamically irrelevant (wandering) domains. This produces an
interesting class of surface homeomorphisms which includes torus
Anosovs and Thurston's pseudo-Anosov maps. These maps --- called {\em
generalized pseudo-Anosov maps} --- preserve a pair of invariant
measured foliations with (possibly infinitely many) singularities,
giving the underlying surface a naturally defined complex
structure. Regarded from the point of view of this complex structure,
the invariant foliations become the horizontal and vertical
trajectories of an integrable quadratic differential (with possibly
infinitely many poles) and the generalized pseudo-Anosov map becomes a
Teichm\"uller mapping. 

The construction of generalized pseudo-Anosov maps is done by first
introducing generalized train tracks: smooth branched 1-submanifolds, with
possibly infinitely many branches, of the
ambient surface. We describe how to find invariant generalized train
tracks for certain surface homeomorphisms which, together with measure
theoretic information obtained from transition matrices, are the main
ingredients in the construction of generalized pseudo-Anosov
maps. This is a variant of the same circle of ideas mentioned above.

The paper is organized as follows: in Section~\ref{sec:natext} the
natural extension is defined, the class of {\em thick graph maps} is
introduced and, in Proposition~\ref{prop:natext}, the circle of ideas
mentioned above is closed. In Section~\ref{sec:0ent} the 0-{\em
entropy} equivalence relation, which collapses dynamically irrelevant
domains, is introduced and a theorem (Theorem~\ref{thm:0ent}) is
stated about the quotient maps so obtained --- this is the first way
in which generalized pseudo-Anosov maps appear. In Section~\ref{sec:gentt}
infinite train tracks are defined as the main tool to construct
generalized pseudo-Anosov maps in the following section. The complex
structure is discussed in the last subsection of Section~\ref{sec:genpA}.

The writing style is Sullivanian as is fitting for a paper prepared for such
an occasion. Many arguments are only sketched and some are omitted entirely.

\medskip
\noindent
{\bf Acknowledgements:} I have worked on or around the main subjects
in this paper for many years. Sevaral people in the course of several
conversations have influenced me. Dennis Sullivan, of course, was a
major influence as my adviser. Others who have directly contributed to
the ideas in this paper are Michael Handel, with whom I had
many conversations about many things, including infinite train tracks
and generalized pseudo-Anosovs; Fred Gardiner, who taught me
Teichm\"uller theory and how to show a puncture is a puncture; and
my friend and coworker Toby Hall, whose work was the original stimulus
for me to study the subject and who wrote a good part of
Section~\ref{sec:gentt} when we were engaged in a project which has
now been put on hold.
 
\section{The natural extension}
\label{sec:natext}

In this paper a {\em dynamical system} will mean a
continuous self-map of a topological space, at least. As we go along,
we may require more of our maps, for example, that they be differentiable
or diffeomorphisms.

\subsection{Definition and examples}

Let $f\colon X\to X$ be a continuous surjective map of a topological
space $X$.
If $f$ is not invertible, there is a naturally
associated invertible map: set 
\[ \hat{X}=\{(x_0,x_1,x_2,\ldots)\in \prod^\infty_0X; f(x_{i+1})=x_i,
\mbox{for}\ i=0,1,2,\ldots\} \] and define
$\hat{f}\colon\hat{X}\to\hat{X}$ setting $\hat
f(x_0,x_1,x_2,\ldots)=(f(x_0),x_0,x_1,\ldots)$. This map is called the
{\em natural extension} of $f$. The space $\hX$ is also known as the
{\em inverse} or {\em projective limit space} and the map $\hat f$ as
the {\em inverse} or {\em projective limit map} associated to $f\colon
X\to X$.

Here are some prototypical examples of natural extensions which come
up in a variety of contexts in the study of low-dimensional dynamical
systems.

\begin{examples}
\begin{enumerate}[a)]
\item Let $X=\mathbb{S}^1$ be the unit circle in the complex plane and 
$f\colon\mathbb{S}^1\to\mathbb{S}^1$ be the squaring map $f(z)=z^2$ in
complex notation. Then $\widehat{\sone}$ is the {\em dyadic
solenoid}. It is a fiber bundle over $\sone$ with fiber a dyadic Cantor
set. 
\item Let $\mathbb{C}^*= \mathbb{C}\setminus 0$ and
$f\colon\mathbb{C}^*\to\mathbb{C}^*$ be again $f(z)=z^2$. Then
$\widehat{\mathbb{C}^*}$ is the {\em complex} dyadic solenoid and is a
fiber bundle over $\mathbb{C}^*$ with same fiber as above. 
\item A variant of the previous example may be obtained by
restricting $f$ to $X=\mathbb{C}\setminus\mathbb{D}$, the exterior of
the closed unit disk. Since the 
action of $f$ on $X$ has a fundamental domain (namely, any annulus of
the form $A=\{z\in\ofc\,;1<R\le |z|\le R^2\}$), so does the action of $\hf$ on $\hX$
and we can take the quotient $\cS=\hX/\hf$. This is {\em
Sullivan's Riemann Surface Lamination}. $\cS$ is a {\em lamination}
as were the previous examples, but, unlike them,
for which collapsing transversals produces a good old space,
this time the quotient is the {\em branched surface} obtained as the
quotient $A/f$ of the annulus by identifying its boundary under the
dynamics.    
\item Let $I=[0,1]$ and $f$ be the {\em tent map} $f(x)=2x$ if $x\in
[0,1/2]$ and $f(x)=2-2x$ if $x\in [1/2,1]$. The inverse limit space
$\hat{I}$ is the {\em Knaster continuum}. Again, $\hat I$ is a
laminated space and the quotient by collapsing transversals is $I$,
but $\hat I$ is not a fiber bundle over $I$, because $0$ and $1$ have
irregular fibers. We will talk more about examples like this below.
\end{enumerate}
\end{examples}

\subsection{The natural extension in 2-dimensional dynamics}

In this section we introduce a class of maps which are suitable for
several applications (see~\cite{BH,dCH1,FM}).  They are called {\em
thick graph maps} and, as the name suggests, they are essentially
graph endomorphisms that have been thickened and made into surface
homeomorphisms.  All of their interesting dynamics is contained in a
subsurface called a {\em thick graph}, that is, a graph in which each
point has been thickened up, either to a disk or an arc according as
the point is a vertex or a regular point of the graph. The
homeomorphisms are assumed to act on the thick graph in such a way
that they induce endomorphisms of the underlying graph.

\begin{figure}[htbp]
\psfrag{1}[t]{$1$-junction}
\psfrag{2}[t]{$2$-junction}
\psfrag{3}[t]{$3$-junction}
\psfrag{4}[b]{leaves}
\begin{center}
\pichere{0.6}{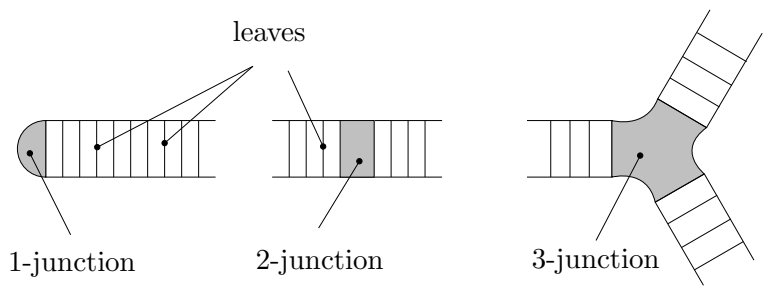}
\caption{Charts in a thick graph.}
\label{fig:chart}
\end{center}
\end{figure}

\begin{defns}
A {\em thick graph} is a pair $(S,\ofg)$, where $S$ is a closed
orientable surface endowed with a fixed metric compatible with its
topology and $\ofg$ is a compact subsurface of $S$ (with boundary) which
is partitioned into compact {\em decomposition elements}, such that
\begin{enumerate}[i)]
\item Each decomposition element of $\ofg$ is either a {\em leaf}
homeomorphic to $[0,1]$, or a {\em junction} homeomorphic to $D^2$
(the unit disk in $\R^2)$.
\item The boundary in $\ofg$ of each junction is a finite number of
disjoint arcs: if there are $k$ such arcs, then the disk is called
a {\em $k$-junction}.
\item The set of $k$-junctions with $k\not =2$ is finite.
\item Each decomposition element which is not in the accumulation of a
sequence of distinct $2$-junctions is contained in a chart as depicted
in Figure~\ref{fig:chart}.
\item Each component of $S\setminus\ofg$ is an open disk.
\end{enumerate}
If $(S,\ofg)$ is a thick graph, let $\sim$ be the equivalence relation
on $\ofg$ given by $x\sim y$ if and only if $x$ and $y$ lie in the
same decomposition element. Then $G=\ofg/\!\!\sim$ is a graph, whose
vertices (which may have valence $2$) correspond to the junctions of
$\ofg$: the canonical projection will be denoted $\pi\colon\ofg\to
G$. The vertex set of $G$ will be denoted $V$, and the union of the
junctions of $\ofg$ will be denoted $\ofv$: thus $\ofv=\pi\I(V)$. The
components of $\ofg\setminus\ofv$ are called {\em strips}: each strip
is therefore homeomorphic to $(0,1)\times[0,1]$.  The union of the
closures of the strips will be denoted by $\ofe$ and the corresponding
set in the quotient, the set of {\em edges} of $G$, will be denoted by
$E$. Thus $\ofe\cap\ofv$ consists of a collection of closed arcs which
are the boundary components (in $ofg$) of both the junctions and the
strips of $\ofg$.
\end{defns}

\begin{rmks}
\begin{enumerate}[a)]
\item We will often refer to the thick graph as $\ofg$ alone,
although it should be kept in mind that there is a surface floating
around. 
\item Notice that, at this point, we allow a thick graph to have
infinitely many 2-junctions and, therefore, infinitely many strips
(and, thus, the graph $G$ to have infinitely many edges).
\end{enumerate}
\end{rmks}

If $(S,\ofg)$ is a thick graph and $F\colon (S,\ofg)\to(S,\ofg)$
is a homeomorphism (i.e., a homeomorphism $F\colon S\to S$ with
$F(\ofg)\subset\ofg$) under which the image of each decomposition
element of $\ofg$ is contained in a decomposition element, then
$F\vert_\ofg$ induces a graph endomorphism $f\colon G\to G$
such that $\pi\circ F\vert_{\ofg} =f\circ\pi$.

\begin{defn}
A {\em thick graph map of $(S,\ofg)$} is an orientation-preserving
homeomorphism $F\colon(S,\ofg)\to(S,\ofg)$ such that:
\begin{enumerate}[i)]
\item $F(\ofg)\subset\Int(\ofg)$.
\item If $\gamma$ is a decomposition element of $\ofg$, then
$F(\gamma)$ is contained in a decomposition element, and
$\diam(F^n(\gamma))\to 0$ as $n\to\infty$. 
\item The induced graph endomorphism $f\colon G\to G$ is 
piecewise monotone (that is, there is a finite subset $L$ of $V$ such
that $f\I(x)\cap U$ is connected for each $x\in G$ and each component
$U$ of $G\setminus L$); and is strictly monotone away from the
preimages of vertices (that is, every $x\in G\setminus f\I(V)$ has a
neighborhood on which $f$ restricts to an embedding).
\item For each component $U$ of $S\setminus\ofg$ there is a (least)
positive integer $n_U$ for which either $F^{n_U}(U)\subset\ofg$ or
$F^{n_U}(U)\cap U\not=\emptyset$, in which case $U$ contains a period
$n_U$ point $p_U$ of $F$, which is a source whose immediate basin
contains $U$: that is, $F^{-kn_U}(x)\to p_U$ as $k\to\infty$ for all
$x\in U$.
\end{enumerate}
\end{defn}

\begin{rmks}
\begin{enumerate}[a)]
\item Item iv) in the definition says that the dynamics of a thick
graph map in $S\setminus\ofg$ is easily understood and
uninteresting. 
\item If $F\colon (S,\ofg)\to (S,\ofg)$ is a thick graph map, then so
are all of its forward iterates $F^n$ ($n\ge1$), and the graph endomorphism
induced by $F^n$ is $f^n\colon G\to G$.
\item Let \gm be the quotient of a thick graph map. A point $x\in G$
is a {\em critical point} if $f$ is not a local homeomorphism
at $x$. Since thick graphs may have infinitely many
2-junctions, the forward orbit of critical points of $f$ may be infinite.

\end{enumerate}
\end{rmks}

\begin{example}\normalfont
\label{ex:hs}
The first example is Smale's horseshoe map which will be denoted by
$F_1\colon(\sph,\ofi)\to(\sph,\ofi)$ here and in what follows. It is
shown in Figure~\ref{fig:hs}. The thick graph in this case is 
a thick tree --- a thick interval, in fact --- and is denoted by $\ofi$. The
point at infinity in $\sph$ is a repeller whose basin contains all points
outside $\ofi$. The horseshoe has two saddle fixed points which are
labelled $x_0$ and $x_1$ (shown as $\bullet$ and $\circ$,
respectively, in Figure~\ref{fig:hs}) and an attracting fixed point in
the 1-junction on the left denoted by $x$ (shown as
\rule{1.5mm}{1.5mm}). The quotient tree is the interval and the
quotient map --- the `flat top' tent map $f_1\colon I\to I$ --- is also
shown in the figure: the image is shown slightly separated so it is
possible to see what $f_1$ does to the interval. This way of
representing graph maps will always be used in what follows.
\begin{figure}[htbp]
\begin{center}
\lab{p}{\pi}{}
\lab{T}{\ofi}{l}
\lab{F}{F_1}{}
\lab{f}{f_1}{}
\lab{x}{x}{}
\lab{y}{x_0}{}
\lab{z}{x_1}{}
\pichere{0.4}{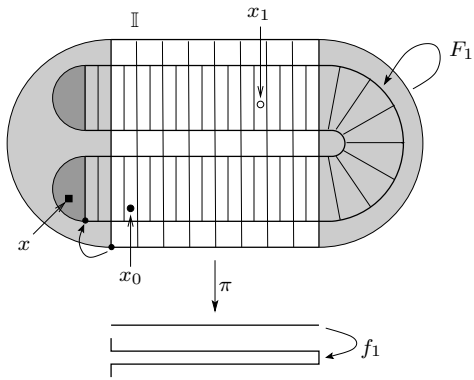}
\caption{The horseshoe map.}
\label{fig:hs}
\end{center}
\end{figure} 
\end{example}

\begin{example}
\label{ex:thtrmap}
In this example the map is again a sphere homeomorphism and the thick graph
is again a thick tree as shown in Figure~\ref{fig:thtrmap}. It is denoted
$F_2\colon(\sph,\oft)\to(\sph,\oft)$. The 1- and
2-junctions of $\oft$ contain a periodic orbit of period 6 ($\circ$) and the 3-junctions
contain a periodic orbit of period 2 ($\bullet$). The quotient tree and tree
endomorphism $f_2\colon T\to T$ are shown in Figure~\ref{fig:trmap}.
\begin{figure}[htbp]
\begin{center}
\lab{T}{\oft}{}
\lab{F}{F_2}{}
\pichere{0.6}{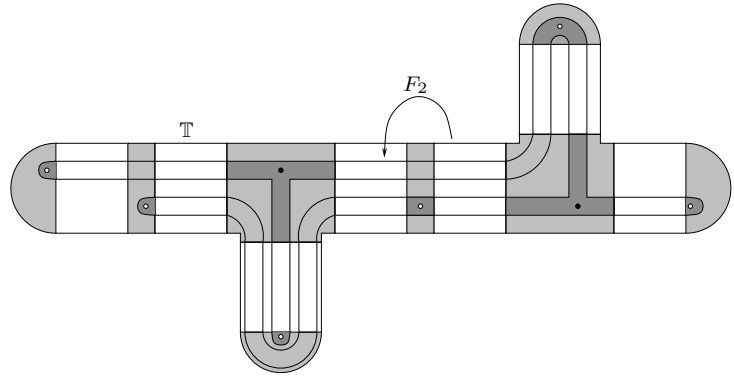}
\caption{A thick tree map.}
\label{fig:thtrmap}
\end{center}
\end{figure}

\begin{figure}[htbp]
\begin{center}
\lab{1}{1}{b}\lab{2}{2}{b}\lab{3}{3}{b}\lab{4}{4}{b}\lab{5}{5}{}\lab{6}{6}{}
\lab{7}{f(1)}{}\lab{8}{f(2)}{}\lab{9}{f(3)}{}\lab{a}{f(4)}{l}\lab{b}{f(5)}{l}
\lab{c}{f(6)}{}  
\lab{f}{f_2}{}
\pichere{0.6}{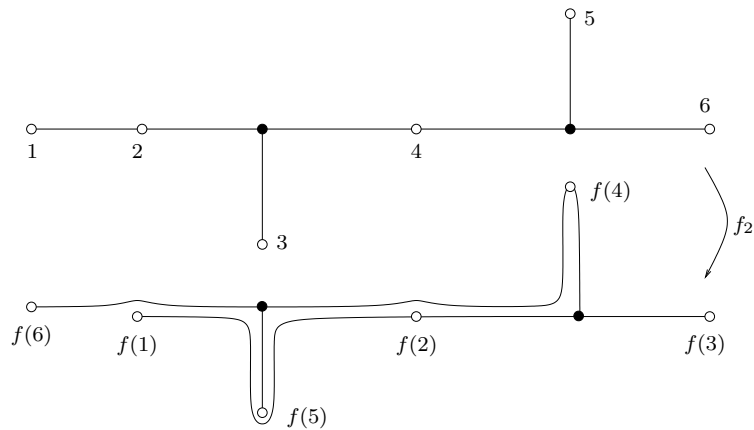}
\caption{The quotient tree map of the thick tree map in
Figure~\ref{fig:thtrmap}.}
\label{fig:trmap}
\end{center}
\end{figure}
\end{example}

\begin{rmk}
Both maps $F_1$ and $F_2$ above can be made to be diffeomorphisms. This
will be used below when we talk about stable and unstable manifolds of
their periodic points.
\end{rmk}

Let $F\colon\ofg\to\ofg$ be a thick graph map and $f\colon G\to G$ its
quotient.  The infinite nested intersection $\Lambda=\bigcap_{n=0}^\infty
F^n(\ofg)$ is a compact subset of $\ofg$ which is invariant under
$F$. We now relate the dynamics of $F$ on $\Lambda$ with the natural
extension of $f$ (see~\cite{Ba}).

\begin{prop}
\label{prop:natext}
The maps $F|_\Lambda\colon\Lambda\to\Lambda$ and $\hf\colon\hat G\to
\hat G$ are topologically conjugate, that is, there exists a
homeomorphism $\hat\pi\colon\Lambda\to\hat G$ such that $\hat\pi\circ
F=\hat G\circ\hat\pi$. 
\end{prop}
\begin{proof}
Let $\pi\colon\ofg\to G$ be the projection.  The map $\hat\pi$ is defined
setting \[\hat\pi(z)=(\pi(z),\pi(F^{-1}(z)),\pi(F^{-2}(z)),\ldots)\] for each
$z\in\Lambda$. It is straightforward to check that $\hat\pi$ is well
defined, continuous and surjective. Injectivity follows from the assumption
that $\diam(F^n(\gamma))\rightarrow 0$ as $n\rightarrow\infty$.
\end{proof}

\section{Another kind of quotient} 
\label{sec:0ent}

In this section we describe a way of modifying two-dimensional maps
by a semi-conjugacy which collapses `irrelevant' dynamics: that is,
parts of the space which do not carry entropy. The semi-conjugacy is
defined quite generally and the space of maps it yields is quite
interesting in its own right. A more thorough treatment of the
equivalence relation can be found
in~\cite{dCP}. The space of quotient homeomorphisms will be treated in
a forthcoming paper. 

\subsection{The 0-entropy equivalence relation}

We start by recalling Bowen's definition of topological
entropy~\cite{Bowen}. If $X$ is a metric space and $F\colon X\to X$ is
a uniformly continuous map, we say that $x,y\in X$ are
$(n,\epsilon)$-{\em separated} if it is possible to distinguish
between the orbits of $x$ and $y$ up to $n-1$ iterates with precision
$\epsilon$. That is, $x,y\in X$ are $(n,\epsilon)$-separated if
$d(F^j(x),F^j(y))>\epsilon$ for some $0\leq j<n$. The topological
entropy of $F$ is defined to be the limit as $\epsilon\to 0$ of the
exponential growth rate of the number of $(n,\epsilon)$-separated
orbits as $n\to\infty$. If $K\subset X$ is a compact subset and we
only count those orbits which start in $K$, we obtain the entropy of
$F$ {\em in} $K$, denoted $h_F(K)$. More precisely, if we denote by
$s(n,\epsilon,K)$ the cardinality of a maximal
$(n,\epsilon)$-separated subset of K, then
\[h_F(K)=\lim_{\epsilon\ra 0}\limsup_{n\ra\infty}\frac{1}{n}
\ln s(n,\epsilon,K)\]
and the entropy of $F$ is defined by
\[h(F)=\sup\{h_F(K); K\subseteq X, K\ \mbox{compact}\}.\]

\begin{defn}
If $F\colon X\to X$ is a homeomorphism, we define two points $x$ and
$y$ to be {\em $0$-entropy equivalent} if there is a continuum (that
is, a compact connected set) $K$ which contains both points and for
which
\[h_F(K)=0=h_{F\I}(K).\]
\end{defn}

\begin{rmks}
\begin{enumerate}[a)]
\item That this indeed defines an equivalence relation follows from
two facts: 1)~$h_F(K\cup K')\le\max\{h_F(K),h_F(K')\}$ and 2)~the
union of two continua containing a point in common is also a
continuum.
\item Notice that if $K$ is a proper subset of $X$, it is not necessarily the
case that $h_F(K)=h_{F\I}(K)$.
\item If $F$ is not invertible we can consider the equivalence
relation defined using only the first equality above. It would be
interesting to understand this equivalence relation --- and the ones
mentioned below --- for interval or, more generally, tree
endomorphisms.
\item In general, we can consider the family of equivalence relations
$\sim_\alpha$, indexed by a positive real $\alpha$,
declaring two points to be $\sim_\alpha$-equivalent if there is a continuum
containing both and carrying entropy strictly smaller than $\alpha$.
\end{enumerate}
\end{rmks}

The 0-entropy equivalence relation is most interesting for two-dimensional
systems. 
\begin{example}
\label{ex:hs0ent}
Let us describe its equivalence classes for the horseshoe
map~ $F_1$.  Denote by ${\mathcal H}^{u}$ and ${\mathcal H}^s$ the closures of
the unstable and stable manifolds of the fixed point~$0$ (or
indeed of any other periodic point, since their closures coincide) and
let ${\mathcal H} = {\mathcal H}^s\cup{\mathcal H}^u$. Equivalence classes are of
four kinds: 
\begin{enumerate}[a)]
\item Closures of connected components of $\sph\setminus{\mathcal H}$.
\item Closures of connected components of ${\mathcal H}^u\setminus{\mathcal H}^s$
(not already contained in sets in a)). 
\item Closures of connected components of ${\mathcal H}^s\setminus{\mathcal H}^u$
(not already contained in sets in a)).
\item Single points which are in none of the sets in a), b) or c).
\end{enumerate}
\end{example}

To see that these sets do not carry entropy, notice that all points in
any connected component of $\sph\setminus{\mathcal H}$ (before taking the
closure) converge to the attracting fixed point $x$. It is not hard to
see that, after taking the closure nothing more exciting happens and
this shows the sets in a) indeed carry no entropy. The same holds for
sets of types b) and c). To see that any larger continuum must contain
entropy, notice that if $C$ is a connected set that contains two
distinct sets among the ones described above, then it must intersect a
Cantor set's worth of invariant manifolds, either stable or unstable
(or both). It follows that one of its $\omega$- or $\alpha$-limit sets
contains all the nonwandering set of the horseshoe and therefore one of 
$h_F(C)$ or $h_{F\I}(C)$ equals $\ln 2$. 

In~\cite{dCP} it is shown that if $F$ is a $C^{1+\epsilon}$ surface
diffeomorphism, then the 0-entropy equivalence classes form an {\em
upper semi-continuous monotone} decomposition of the surface. In particular,
the quotient space is a {\em cactoidal} surface (roughly, a surface
with nodes; see~\cite{Moo,RoSt}).
Since the equivalence is dynamically defined, $F$ projects to
a homeomorphism $F/\!\!\sim$ on the quotient space. Moreover, any
nontrivial continuum in the quotient space carries
entropy of either the quotient homeomorphism or its inverse. The
quotient map by the 0-entropy equivalence relation should be thought
of as a `tight' version of the original map in which all the wandering
domains have been collapsed to points.  

\begin{figure}
\pichere{0.5}{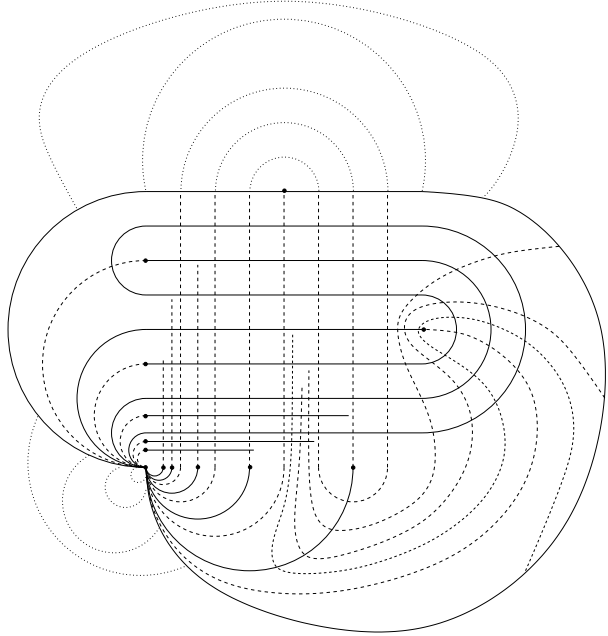}
\caption{The quotient of the sphere under the 0-entropy equivalence
relation for the horseshoe.}
\label{fig:tightHS}
\end{figure}

The quotient of the sphere by
the $0$-entropy horseshoe equivalence of Example~\ref{ex:hs0ent} is
shown in Figure~\ref{fig:tightHS}. The quotient space is a sphere
(we are collapsing everything outside the homoclinic tangle to a
point), obtained by identifying the solid boundary in the figure along
the dotted arcs from the mid-point at the top to the
corner point on the lower left. The stable and unstable manifolds of
the horseshoe project to two transverse foliations with singularities,
represented by solid and dashed lines, respectively. In fact, these
foliations carry transverse invariant measures whose product gives a
measure on the sphere. The quotient map preserves both foliations,
dividing one of the transverse measures by 2 and multiplying the other
by 2, so that the product measure is invariant. This map is a
{\em generalized pseudo-Anosov map}, which will be defined presently.
We will also state a theorem that generalizes this construction to a
class of thick graph maps.

\subsection{The Markov story and generalized pseudo-Anosov maps}
\label{subsec:markov}

We start by introducing the basic concepts in the Perron-Frobenius theory
for non-negative matrices (see~\cite{Baladi,Gant}).
Let $M$ be a square matrix with non-negative integer entries. $M$ is
said to be {\em reducible} if, by a permutation of the index set, it
is possible to put it in triangular block form:
\[ M = \left[ \begin{array}{cc}
               A & 0 \\
               B & C
             \end{array} \right] . \]
Otherwise, $M$ is said to be {\em irreducible}. The matrix $M$ is said
to be {\em irreducible and aperiodic} if there exists a positive
integer $k$ such that $M^k$ is positive, that is, all its entries are
positive.  A non-negative irreducible matrix $M$ has a unique positive
eigenvector (up to scaling) and the associated eigenvalue $\lambda$
--- called the {\em Perron-Frobenius eigenvalue} of $M$ --- equals the
spectral radius of $M$. If $\lambda=1$ then $M$ is a cyclic
permutation matrix. Otherwise, $\lambda>1$ and, for every $i,j$, there
exists a power $k$ such that the $ij$-entry of $M^k$ is arbitrarily
large (in fact, the entries of $M$ grow like (const)$\times\lambda^k$). If 
$M$ is irreducible then $\lambda$ is a simple root of the
characteristic polynomial of $M$ and if $M$ is also aperiodic, then
$\lambda$ is the only eigenvalue on the circle $\{z\in \mathbb{C}\,;
|z|=\lambda\}$.  

\begin{defns}
A thick graph map $F\colon\ofg\to\ofg$ with quotient $f\colon G\to G$
will be called {\em Markov} if it satisfies the following additional
conditions:
\begin{enumerate}[i)]
\item $\ofg$ has a finite number of strips and junctions (i.e., the
graph G has finitely many edges). For each strip $s$ we
fix a homeomorphism $h_s \colon \sb \to [0,1]\times[0,1]$ from the closure
of the strip to the closed unit square, so that the decomposition
elements of $s$ are of the form ${h_s}^{-1}(\{x\}\times [0,1])$, for
$0<x<1$.
\item $F$ is linear with respect to the structure homeomorphisms $h_s$,
that is, in each connected component of $s_i \cap F^{-1}(s_j)$, where
$s_i$, $s_j$ are strips, $F$ contracts vertical coordinates uniformly
by a factor $\mu_{ij}<1$ and expands horizontal coordinates uniformly
by a factor ${\lambda}_{ij} \ge 1$.
\item If $J, J'$ are junctions such that $F(J) \subset J'$
then $F(\partial_\ofg J) \subset \partial_\ofg J'$. 
\item If $J$ is a periodic junction of least period $n$, then $J$ has
an attracting periodic point of least period $n$ in its interior whose
basin contains $\Int(J)$.
\end{enumerate}
We can associate a {\em transition matrix} $M=[m_{ij}]$ to a Markov thick
graph map: letting $E=\{e_1,e_2,\ldots,e_n\}$ be the edges of $G$, set 
\[ m_{ij}=\mbox{number of times\ } f(e_j) \mbox{\ crosses\ }
e_i.\] 
\end{defns}


\begin{rmks}
\begin{enumerate}[a)]
\item Notice that Markov thick graph maps can be made differentiable
and we will assume, whenever we talk about them, that they are
diffeomorphisms of the surface $S$. In particular, we will talk freely
about stable and unstable manifolds of their periodic points.
\item Let $F\colon (S,\ofg)\to(S,\ofg)$ be
a Markov thick graph map whose transition matrix is irreducible and
aperiodic, let $\Lambda=\bigcap_{n=0}^\infty F^n(\ofg)$ and
$p\in\ofg$ be any saddle periodic point of $F$. Then it is easy to see that
$\Lambda=\ol{W^u(p)}$.
\end{enumerate}
\end{rmks}

\begin{figure}[htbp]
\begin{center}
\pichere{0.6}{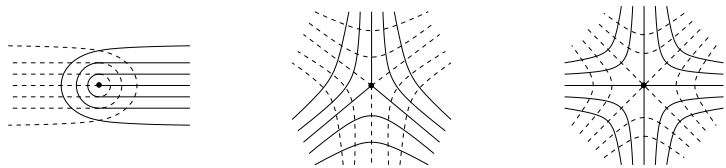}
\caption{Pronged singularities of the invariant foliations.}
\label{fig:prongs}
\end{center}
\end{figure}

\begin{defn}
A surface homeomorphism $\Phi\colon S\to S$ is called a {\em generalized
pseudo-Anosov map} if it satisfies the following
conditions: there exist a pair $(\cF^s,\mu^s)$, $(\cF^u,\mu^u)$ of
transverse measured foliations with singularities --- either modeled on
pronged singularities as in Figure~\ref{fig:prongs} or
accumulations of such, of which there are only finitely many --- and a
real number $\lambda>1$ such that
\[\begin{array}{ccc}
\Phi(\cF^s,\mu^s) &=& (\cF^s,\lambda\mu^s) \\
\Phi(\cF^u,\mu^u) &=& (\cF^u,\frac{1}{\lambda}\mu^u). 
\end{array}\]
\end{defn}

Examples of generalized pseudo-Anosov maps include the torus Anosov
maps and Thurston's pseudo-Anosov maps. We argue below that Markov
thick graph maps also give rise to generalized pseudo-Anosovs. The
definition above, however, probably includes many other maps and these
Markov examples should form a dense set in the space of all
generalized pseudo-Anosovs. A more thorough study of these issues,
including a description of a uniform structure on the set of
generalized pseudo-Anosovs, is currently under way and will hopefully
appear in forthcoming papers.

\begin{thm}
\label{thm:0ent}
Let $F\colon(S,\ofg)\to(S,\ofg)$ be a Markov thick graph map with
$\ofg$ of the homotopy type of $S$ minus a point.
Assume the associated transition matrix is irreducible and aperiodic.
Then the quotient of $S$ by the 0-entropy equivalence relation is
homeomorphic to $S$ and $F$ projects to a generalized pseudo-Anosov
homeomorphism $\Phi\colon S\to S$. 
\end{thm}

A proof of this theorem shall appear in a forthcoming paper.  A more
detailed account of the 0-entropy equivalence relation and its
quotients is given in~\cite{dCP}. Below we will present an alternative
construction of the generalized pseudo-Anosov quotient maps.

\section{Generalized train tracks}
\label{sec:gentt}

In this section we will only deal with {\em finite} thick graphs, that
is, thick graphs with finitely many strips and junctions. To simplify
the exposition and the statements, we also assume that if $(S,\ofg)$
is a thick graph then $\ofg$ has the homotopy type of the punctured
surface $S\setminus\{p\}$, where $p\in S\setminus\ofg$. We refer to
$p$ as the {\em point at infinity} and denote it by $\infty$. In
general, $\ofg$ has the homotopy type of the several times punctured
surface $S\setminus\{p_1,\ldots,p_k\}$. If this is the case, the
discussion below has to be appropriately modified, but the ideas are
essentially the same. The statements made here hold for the once
punctured case.

\subsection{Definitions}

We start by fixing some notation and presenting some definitions. 
Suppose that $f\colon X\to X$ is a homeomorphism. Then $f$ is said
to be {\em supported} on a subset $U$ of $X$ if $f$ is the identity on
$X\setminus U$. A second homeomorphism $g\colon X\to X$ is {\em
isotopic to $f$} if there is a continuous map $\psi\colon
X\times[0,1]\to X$ such that each {\em slice map} $\psi_t\colon X\to
X$ defined by $x\mapsto \psi(x,t)$ is a homeomorphism, and $\psi_0=f$
and $\psi_1=g$. The map $\psi$ is called an {\em isotopy} from $f$ to
$g$.  A {\em pseudo-isotopy} is a continuous map $\psi\colon X\times
[0,1] \to X$ such that, for $0 \le t < 1$, the slice maps $\psi_t$ are
homeomorphisms onto their images.  The isotopy or pseudo-isotopy is
said to be {\em supported} on a subset $U$ of $X$, denoted
$\mbox{supp}(\psi)=U$, if the homeomorphisms $\psi_t$ are all equal on
$X\setminus U$, and is said to be {\em relative} to $U$ if
$\psi_t(U)\subset\psi_0(U)$ for all $t\in[0,1]$. A map $g\colon X\to
X$ is called a {\em near-homeomorphism} if it can be arbitrarily well
approximated by homeomorphisms. Thus, if $\psi$ is a pseudo-isotopy,
the map $g=\psi_1$ is a near-homeomorphism.

Let $(S,\ofg)$ be a finite thick graph, and $A\subset\ofg$ be a finite
set, each of whose points lies in the interior of a junction and no
junction contains more than one point of $A$. For each strip $s$ of
$\ofg$, let $\gamma_s$ be an arc joining the two boundary components
of $s$ in $\ofg$ and intersecting each leaf $s$ exactly once. In terms
of the structure homeomorphisms $h_s$ introduced in
Subsection~\ref{subsec:markov}, we can take $\gamma_s =
h_s^{-1}([0,1]\times\{1/2\})$. Let $R \subseteq\ofg$ be the union of
the arcs $\gamma_s$. The endpoints of the arcs $\gamma_s$ are called
{\em switches} and we denote by $L$ the set of switches.

\begin{defns}
Let $\tau\subseteq\ofg\setminus A$ be a graph with vertex set~$L$ and
countably many edges, each of which intersects $\partial\ofv$ only
at~$L$, such that
\begin{enumerate}[i)]
\item $\tau\cap\ofe = R$, and
\item No two edges $e_1$, $e_2$ contained in a given junction $J$ are
{\em parallel}: that is, they do not bound a disk which contains no
point of $A$ or other edges.
\end{enumerate}

The isotopy class of $\tau$ by isotopies supported on $\ofv\setminus
A$ (the set of junctions of $\ofg$ minus the points in $A$) is called
a {\em generalized train track}\footnote{To be painfully precise, we
should talk about the isotopy class of the inclusion map
$\iota\colon\tau\hookrightarrow\ofg$, but we won't do it.} for
$(\ofg,A)$. We will always refer to $\tau$ itself as the generalized
train track, but it should be kept in mind that we do not
distinguished between $\tau$ and $\tau'$ if it is possible to deform
one to the other without crossing over points of $A$.

The edges of $\tau$ which are contained in $\ofe$ (that is, the
connected components of $R$) are called {\em
real}, and those which are contained in $\ofv$ are called {\em
infinitesimal}. Write $I$ for the set of infinitesimal edges of
$\tau$.

A generalized train track~$\tau$ is {\em finite} if it has only
finitely many edges. An infinitesimal edge is called a {\em bubble} if
its two endpoints coincide. An edge of $\tau$ is {\em homotopically
trivial} if it is a bubble which bounds a disk containing no point of
$A$, and is {\em homotopically non-trivial} otherwise.

Clearly a generalized train track $\tau$ for $(\ofg,A)$ is determined
by its infinitesimal edges. It will sometimes be convenient to write
$\tau(I)$ for the generalized train track whose set of infinitesimal
edges is $I$, provided the thick graph and the set $A$ are clear from
the context.

A homotopy of a path $\alpha\colon[0,1]\to X$ is said to be {\em
relative to} $U\subset X$ if the points of $\alpha([0,1])$ that belong
to $U$ do not leave $U$ throughout the homotopy.  Let $\alpha$ be a
homotopy class of paths in $S\setminus A$ relative to
$\partial\ofv$, with endpoints in $\partial\ofv$. Then $\alpha$ is
{\em carried} by a generalized train track $\tau$ if it can be
realized by an edge-path in $\tau$ with alternating real and
infinitesimal edges.
\end{defns}

\begin{rmk}
Although for our purposes a train track is a combinatorial object, we
think of it as a smooth branched 1-submanifold of S. From this
standpoint, the homotopy class of a path is carried by a generalized
train track if there is a smooth representative in the class which is
contained in the train track.
\end{rmk}

Now let $F\colon(S,\ofg,A)\to(S,\ofg,A)$ be a thick graph map such
that $F(A)=A$, where $A=A_F$ is the set of attracting
periodic orbits of $F$.  On each strip $s$ of $\ofg$ define the
pseudo-isotopy $\psi_s\colon \bs
\times [0,1] \to \bs$ to be given in coordinates by $ \psi_s(x,y,t) =
(x, (1-t)y+t/2)$ so that $\psi_s(\cdot,1)$ maps $\bs$ onto
$\gamma_s$. Extend these pseudo-isotopies to a pseudo-isotopy
$\psi_0\colon S\times[0,1]\to S$ in the following way: first extend
the $\psi_s$ to mutually disjoint disk neighborhoods $U_s\supset\bs$,
with $U_s\subset
S\setminus A$, so that they are isotopies
on $U_s\setminus\bs$ and the identity on $\partial U_s$; then extend
them to be the identity elsewhere.

Notice that, if $\tau$ is a generalized train track, then $\psi_0(F(\tau))$
satisfies the definition of train track, except possibly b),
that is, it may have parallel edges. 

\begin{defns}
Define $F_*(\tau)$ to be the generalized train track consisting of
a maximal subset of the edges of $\psi_0(F(\tau))$ which contains
no pair of parallel edges.

The train track $\tau=\tau(I)$ is said to be $F$-{\em invariant} if
$F_*(\tau)$ is isotopic to $\tau$ in $S\setminus A$ and $I$ is minimal
(under inclusion) with this property.
\end{defns}

\begin{rmk}
Notice that $F_*(\tau)$ carries the homotopy class of $F(e)$ for
each edge $e$ of $\tau$.
\end{rmk}

From the definition it follows that there exists a
pseudo-isotopy $\psi$ with the property that $\psi(F(\tau),1)=\tau$.

\medskip
\begin{defn}
The {\em train track map} $\phi\colon\tau\to\tau$ associated to
$F\colon\ofg\to\ofg$ is defined by $\phi(\cdot)=\psi(F(\cdot),1)$.
\end{defn}

\subsection{Construction of invariant generalized train tracks}

Let $F\colon(S,\ofg,A)\to(S,\ofg,A)$ be as before. The following
procedure constructs an $F$-invariant generalized train track $\tau$.

Let $\tau_0=R$, and for each $n\ge 0$ define
$\tau_{n+1}=F_*(\tau_n)$. Since $\tau_0$ is a subset of $\tau_1$,
each $\tau_n$ is naturally isotopic to a subset of $\tau_{n+1}$: hence
each $\tau_{n+1}$ can be adjusted as it is constructed by an isotopy
relative to $A\cup \partial\ofv$ such that $\tau_n\subseteq\tau_{n+1}$. Define
$\tau=\bigcup_{n\ge 0}\tau_n$. Then $\tau$ is a generalized train
track by construction. It easy to see that it is
$F$-invariant.

\begin{rmk}
It is clear that this construction provides the minimal $F$-invariant
generalized train track~$\tau$. Also, given any generalized sub-train track
$\tau'$ of $\tau$, it is clear that the same construction, starting
with $\tau_0=\tau'$, must generate $\tau=\bigcup_{n=0}^\infty
F_*^n(\tau')$. 
\end{rmk}

\begin{example}
\label{ex:hstrtr}
The invariant train track for the horseshoe map $F_1$ is shown in
Figure~\ref{fig:hstrtr}. The set $A$ consists of the fixed point $x$
contained in the left 1-junction of $\ofi$. No bubble encloses it and
it is not shown in the figure.
\begin{figure}
\lab{t}{\tau}{}
\lab{f}{\phi}{b}
\pichere{0.4}{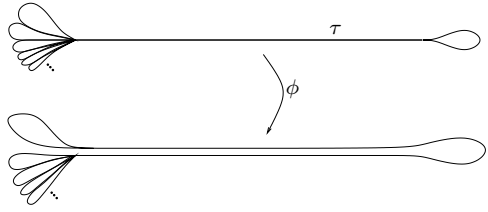}
\caption{The invariant train track for the horseshoe.}
\label{fig:hstrtr}
\end{figure}
\end{example}

\begin{example}
\label{ex:thtrtrtr}
The invariant train track for the thick tree map $F_2$ of
Example~\ref{ex:thtrmap} is shown in
Figures~\ref{fig:thtrtrtr} and ~\ref{fig:bubble}. Here, $A$ consists
of points of both period 6 and period 2 periodic orbits contained in the
junctions of $\oft$. 
\begin{figure}
\lab{t}{\tau}{}
\lab{f}{\phi}{b}
\pichere{0.6}{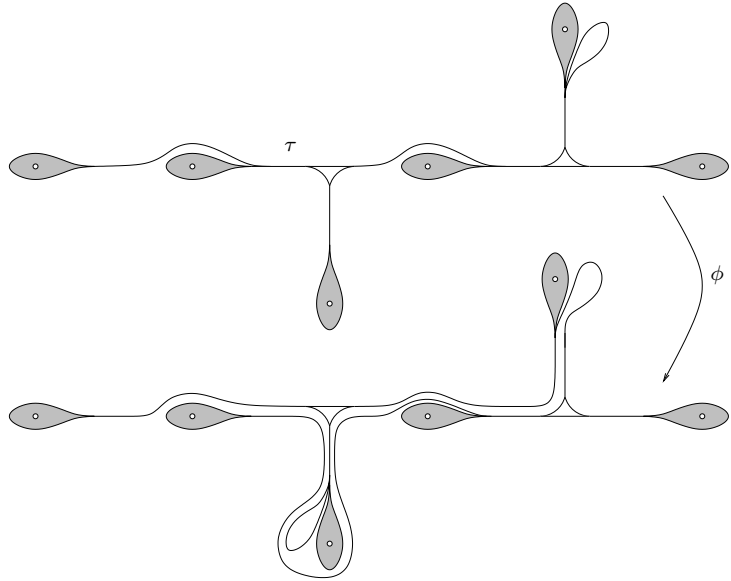}
\caption{The invariant train track for the thick tree map $F_2$ of
Example~\ref{ex:thtrmap}. Figure~\ref{fig:bubble} shows a blow up of
the shaded bubbles in this figure.}
\label{fig:thtrtrtr}
\end{figure}

\begin{figure}
\pichere{0.4}{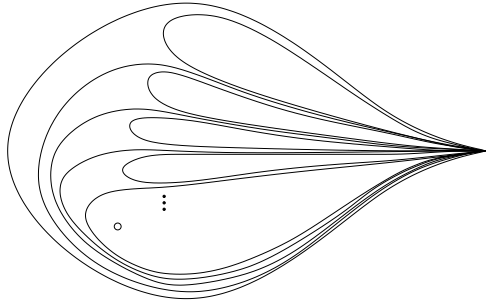}
\caption{Detailed view of the bubbles in
Figure~\ref{fig:thtrtrtr}. There and here $\circ$ represents a periodic point in the
period 6 attracting orbit in Example~\ref{ex:thtrmap}.}
\label{fig:bubble}
\end{figure}
\end{example}

The invariant train track $\tau$ and the train track map $\phi$ should
be thought of as more careful 1-dimensional representations of the
thick graph $\ofg$ and the thick graph map \tgmx: whereas \gm does not
pay attention to junctions --- they are collapsed to points ---
the map $\phi\colon\tau\to\tau$ gives a careful account of the behavior of
the images of strips under iterates of $F$ inside the junctions.

We now gather some useful facts that are straightforward consequences of the
constructions above.
\begin{rmks}
\label{rmks:infinitett}
\begin{enumerate}[a)]
\item $\phi$ is a near-homeomorphism, that is, it can be approximated
arbitrarily well by homeomorphisms. 
\item $\phi$ maps $I$ into itself and each infinitesimal edge
is mapped homeomorphically onto another infinitesimal edge. 
\item At most finitely many
infinitesimal edges are mapped onto a given infinitesimal edge under
$\phi$. 
\item If $e$ is a real edge, then $\phi(e)$ can
only intersect finitely many infinitesimal edges.
\end{enumerate}
\end{rmks}

\begin{figure}[htbp]
\begin{center}
\pichere{0.6}{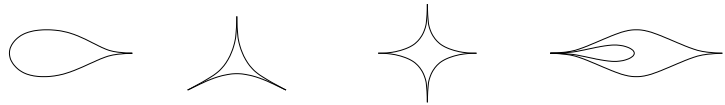}
\caption{Examples of $n$-gons for $n=1,3,4$ and a 1- and a 3-gon together.}
\label{fig:ngons}
\end{center}
\end{figure}

\begin{defn}
By an {\em infinitesimal polygon} we mean a component of the
complement of $\tau$ bounded by finitely many infinitesimal edges. It
is called an {\em $n$-gon} if it is bounded by $n$ infinitesimal
edges (see Figure~\ref{fig:ngons}). 
\end{defn}

\begin{rmks}
\begin{enumerate}[a)]
\item  Bigons (2-gons) are not allowed unless they contain a point of
$A$.
\item It is allowed that an $n$-gon has fewer than $n$ vertices, as
shown in the last diagram in Figure~\ref{fig:ngons}: two vertices in
the 3-gon coincide. 
\end{enumerate}
\end{rmks}

The following proposition is an immediate consequence of the fact that
$\phi$ is a near-homeomor\-phism on a surface. 
\begin{prop}
\label{prop:n-gons}
For each integer $n\ge 1$, $\phi$ induces a 1-1 map on the
collection of $n$-gons. 
If an $n$-gon contains a periodic point in its interior, then it is
periodic under $\phi$ (with the same period). Otherwise, it belongs to
a semi-infinite orbit of $n$-gons $\{\phi^n(\Delta); n\ge 0\}$, where
$\Delta$ is an $n$-gon which is not the image under $\phi$ of any
other $n$-gon. Two such orbits either coincide or are disjoint.
Moreover, in this case, for all but finitely many $n$, $\phi^n(\Delta)$
has at most two vertices.
\end{prop}

\section{Generalized pseudo-Anosov maps}
\label{sec:genpA}

In this section we describe the construction of a generalized
pseudo-Anosov map using an invariant train track, as defined in the
previous section. The construction follows closely the one
in~\cite{BH} for finite invariant train tracks and reduces to that in
case the train track is finite. We also describe an associated complex
structure on the surface with respect to which the map becomes a
Teichm\"uller mapping. We continue to assume that thick graphs have
the homotopy type of a once punctured surface. 

\subsection{The construction}
\label{subsec:construction}

Let \tgm be a thick graph map and $\phi\colon\tau\to\tau$ its
associated invariant train track map as defined in
Section~\ref{sec:gentt}. As before, $R$ and $I$ denote the real
and infinitesimal edges of $\tau$. We number the real edges $\{e_i;
1\le i\le n\}$ and the infinitesimal edges $\{e_i; i>n\}$ and define a
(possibly infinite) transition matrix $M=[m_{ij}]$ setting, as above,
\[ m_{ij}=\mbox{number of times\ } \phi(e_j) \mbox{\ crosses\ } e_i\]

Since $\phi(I)\subset I$, $M$ has block form
\[M=\left[ \begin{array}{cc}
               N_{n\times n} & 0 \\
               B             & \Pi
             \end{array} \right] \]
The matrix $N$ records transitions between real edges and
$\Pi$ records transitions between infinitesimal edges, whereas $B$
records transitions from real to infinitesimal edges.
The (possibly infinite) square matrix $\Pi$ has only 0's and 1's in its
entries. Each of its columns has exactly one non-zero entry and each row has
at most finitely many non-zero entries.
The matrix $B$ (which has $n$ columns and possibly infinitely
many rows) has at most finitely many non-zero entries in each column. 
These observations follow from those in Remarks~\ref{rmks:infinitett}.

\begin{example}
For the horseshoe map, the transition matrix is
infinite but is quite simple: $m_{11}=2$, $m_{i\, i-1}=1$ and
$m_{ij}=0$ otherwise.
\end{example}

\bigskip
\par\noindent{\em Standing Assumption:} It
is assumed throughout this section that the matrix $N$ is irreducible
and aperiodic. This implies that its Perron-Frobenius eigenvalue
$\lambda>1$. It also follows that there is a positive integer $k$ such
that, for every real edge $e$, $\phi(e)$ contains all other real
edges. In particular, $\tau$ is connected. This assumption is not
necessary for all the results that follow, but it simplifies the
discussion. 
\bigskip

We think of $M$ as an operator acting on the space $l^1$ of summable sequences
of real numbers with norm $|y|_1=\sum_{i\ge 1}|y_i|$. From the remarks
above, it follows that $M$ is a bounded operator with 
$\|M\|_1\le\mbox{max}_j\{\sum_i|m_{ij}|\}<\infty$.

Let $Y$ be an eigenvector of $N$ associated to $\lambda$ (and
therefore unique, up to scale). It is
possible to complete $Y$ to an 
eigenvector $y=[Y\: Y']$ of M. Since the columns of $\Pi$ have at most
one non-zero entry which is 1, $\|\Pi\|_1\le 1$ and thus
$\lambda I-\Pi$ is invertible. Setting $Y'=(\lambda I-\Pi)^{-1}BY$ we
have
\begin{eqnarray*}
M \left[ \begin{array}{c}
          Y \\ Y'
         \end{array} \right] &=& \left[ \begin{array}{cc}
                                         N & 0   \\
                                         B & \Pi
                                        \end{array} \right] 
                                \left[  \begin{array}{c}
                                            Y \\ Y'
                                        \end{array} \right] \\
                             &=&  \left[ \begin{array}{c}
                                            NY \\ BY+\Pi Y'
                                         \end{array} \right] \\
                             &=&  \left[ \begin{array}{c}
                                     \lambda Y \\ \lambda Y'
                                         \end{array} \right]
\end{eqnarray*}

In order to see that $Y'$ is a positive vector, notice that 
\begin{eqnarray*}
Y' &=& \frac{1}{\lambda}(I+
\frac{1}{\lambda}\Pi+\frac{1}{\lambda^2}\Pi^2+\ldots)BY  \\
   &=&
\frac{1}{\lambda}(B+\frac{1}{\lambda}\Pi B+\frac{1}{\lambda^2}\Pi^2B+\ldots)Y
\end{eqnarray*}
The matrices $\Pi^kB$ that appear above represent transitions from a
real edge to an infinitesimal edge under the $(k\!+\!1)$-st iterate of the
map $\phi$. In fact, they represent exactly those transitions 
which occurred from a real to an infinitesimal edge in the first
iterate and which then remained among infinitesimal edges
for the next $k$ iterates (the other ways to get from a real to an
infinitesimal edge under the $(k\!+\!1)$-st iterate of $\phi$ are
represented by matrices of the form $\Pi^{k-j-1}BN^j$). But every
infinitesimal edge is the image, under some iterate of $\phi$, of an 
infinitesimal edge of $\tau_1$ and these are the 
intersection of $\phi(\tau_0)$ with $\ofv$, that is, they are
transitions from real
to infinitesimal edges under the first iterate of $\phi$. This
means that, for every $i\ge 1$, there exists $k\ge 1$ such that the
$i$-th row of $\Pi^kB$ is non-zero. Since $Y$ is a positive vector,
it follows that every entry of $Y'$ is non-zero. 

\begin{defn}
A collection $\{y_i=y(e_i)\}_{i\ge 1}$ of non-negative real numbers, called
{\em weights}, is said to satisfy the {\em switch conditions} if, for
each switch $q$ of $\tau$, we have
\[ y(e_{i_0})= \sum w(e_i) + 2\sum w(e_j) \]
where $e_{i_0}$ is the real edge with endpoint at $q$ and the first and
second sums range over the set of infinitesimal edges having one or
both endpoints at $q$ respectively.
\end{defn}

\begin{lemma} 
Let $M$ be the transition matrix associated to
$\phi\colon\tau\to\tau$, $\lambda$ its Perron-Frobenius eigenvalue and
$y=[Y\: Y']=[y_1 y_2 \ldots]$ an eigenvector associated to $\lambda$
as constructed above. Then the set of weights $\{y_i\}_{i\ge 1}$
satisfy the switch conditions.
\end{lemma}

\medskip
\noindent
{\em Proof} (cf. \cite{BH}). 
Fix a large positive integer $k$. The equality $M^ky={\lambda}^ky$
written in coordinates states that for each $i\ge 1$ 
\[ y_i=\frac{1}{{\lambda}^k}\sum_{j}y_j\cdot(\mbox{number of times
${\phi}^k(e_j)$ intersects $e_i$}) \]
If ${\phi}^k(e_j)$ crosses a switch, it must cross edges on both
sides except at its endpoints. Thus, the contribution to both sides of
the switch is the same up to a bounded amount. Letting $k\to\infty$
yields the result. \qed
\medskip

Now, let $X$ be an eigenvector of $N^T$ associated to the
Perron-Frobenius eigenvalue $\lambda$. If we try to complete $X$ to an
eigenvector $[X\: X']$ for the adjoint $M^*$, we are forced to set 
$X'=0$. The reason is that $X'$ is a solution to the equation 
${\Pi}^*X'=\lambda X' \Leftrightarrow (\lambda I-{\Pi}^*)X'=0$. Since 
$\|{\Pi}^*\|_{\infty}\le 1$, the only solution is $X'=0$. 
This is why infinitesimal edges are called such.

We now give a description of the construction of the generalized
pseudo-Anosov homeomorphisms corresponding to \trmx. As was mentioned
before, it follows closely that presented in~\cite{BH}.

Let $x=(x_1,x_2,\ldots,x_n,0,0,\ldots)$ and $y=(y_1,y_2,\ldots)$ be
the eigenvectors of $M^*$ and $M$, respectively, associated to the
Perron-Frobenius eigenvalue $\lambda$ as just described.  To each real
edge $e_i$, $1\le i\le n$ of $\tau$ we associate a Euclidean
rectangle $R_i$ of dimensions $x_i\times y_i$ endowed with foliations
by horizontal and vertical line segments. Each foliation has a
transverse measure induced by Lebesgue measure. The horizontal and
vertical foliations will be called {\em unstable} and {\em stable}
respectively and, under the map, unstable leaves will be stretched and
stable leaves will be contracted by the factors $\lambda$ and
$1/\lambda$, respectively. Place (homeomorphic copies of) these
rectangles on $S$ along the real edges of $\tau$. The infinitesimal
edges of $\tau$ are used to define an equivalence relation on the
vertical sides of the rectangles, as follows. Let $e_j$ be an
infinitesimal edge and $y_j$ the corresponding entry of the
eigenvector $y$. We identify segments of length $y_i$ along the
vertical sides of the rectangles which contain the endpoints of $e_j$
(note that these rectangles could be the same). The facts that the
train track is a subset of $S$ and that the switch conditions are
satisfied imply that there is exactly one way in which these
identifications can be made without self-intersections. The quotient
space $\mathcal R$ of the rectangles under these identifications is a
compact topological surface (with boundary) of the homotopy type of
$S\setminus \{p\}$. It is foliated by unstable `horizontal' lines,
most of which are now infinite, that is, homeomorphic to $\mathbb R$
(the exceptions being those that contain singularities of the
foliations). Moreover, there is defined a map $\tPhi\colon\cR\to\cR$
which stretches the unstable foliation by factor $\lambda$, contracts
the stable foliation (whose leaves are still finite segments) by the
by factor $1/\lambda$ and places $\tPhi(\cR)$ inside $\cR$ in the
manner dictated by \trmx. Restricted to the interior of $\cR$, $\tPhi$
is a homeomorphism, but not along the boundary of $\cR$. Notice that,
by collapsing to points the stable segments, we obtain the graph
$G$. Denoting the projection by $\tilde{\pi}\colon\cR\to G$, $\tPhi$
factors down to
\gmx, that is, $f\circ\tilde{\pi}=\tilde{\pi}\circ\tPhi$. The boundary
$\partial\cR$ has the structure of a smooth finite sided polygon, each
side of which contains a periodic point. We identify segments of
adjacent sides which are mapped to the same segment under some iterate
of \trtrmx. This usually leads to infinitely many identifications of
smaller and smaller pieces of $\partial\cR$ (for example, this is the
case with the horseshoe, as will be seen below).  The points of the
periodic orbit on $\partial\cR$ are all identified and become {\em the
point at infinity.}  Under these further identifications, the quotient
of $\cR$ is homeomorphic to $S$, the `vertical' segments become a
foliation of $S$ by stable leaves (most of which are) homeomorphic to
$\ofr$ and the induced map, denoted by $\Phi\colon S\to S$, now
becomes a homeomorphism: this is the generalized pseudo-Anosov
homeomorphism. The stable and unstable foliations are denoted by
$\cF^s$ and $\cF^u$ respectively. Both are preserved by $\Phi$, the
leaves of the stable foliation being contracted and those of the
unstable foliation being stretched by the factors $1/\lambda$ and
$\lambda$, respectively.


\begin{rmks}
\begin{enumerate}[a)]
\item{The foliations $\cF^s$ and $\cF^u$
are foliations with singularities: to each $n$-gon of $\tau$ there
corresponds an $n$-pronged singularity of the foliations. There may be
infinitely many such singularities, but they can accumulate on at most
finitely many periodic orbits, which are in 1-1 correspondence with a
subset of $A$. In the next section the orbits of 1-pronged
singularities will be studied in greater detail.}
\item The periodic orbits on $\partial\cR$ described above
have been found in different guises by other authors. See, for
example,~\cite{BL,NP}.
\end{enumerate}
\end{rmks}

\begin{example}
The invariant train track for the horseshoe has only one real edge and
therefore the surface ($\sph$) will be constructed by identifying the
edges of one rectangle. The identifications are shown in
Figure~\ref{fig:hsident} by dotted lines. Notice that the equivalence class of the
lower left corner contains infinitely many points. The points marked
with $\bullet$ lie on one orbit of 1-pronged singularities which is
forward and backward asymptotic to the lower left corner, which is
also the point at infinity. The quotient sphere with the two invariant
foliations is shown in Figure~\ref{fig:tightHS}.
\begin{figure}
\pichere{0.5}{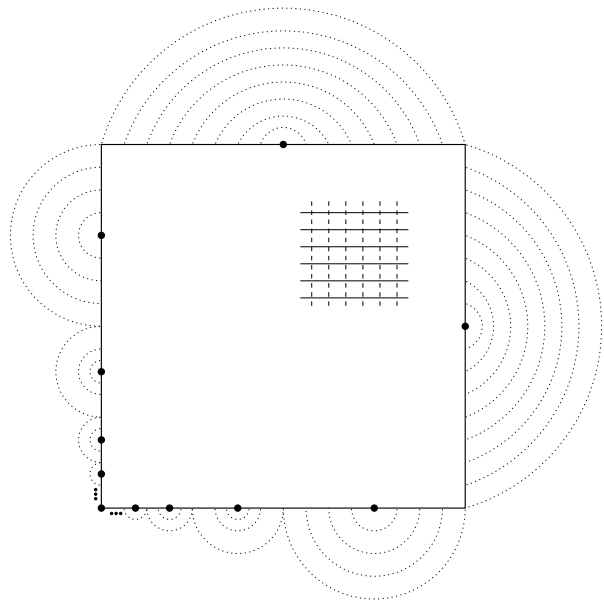}
\caption{The identifications for the horseshoe.}
\label{fig:hsident}
\end{figure}
\end{example}

\begin{example}
In Figure~\ref{fig:thtrident} are shown the identifications on seven
rectangles dictated by the invariant train track for $F_2$. All the
unstable identifications are indicated (dotted lines) whereas only the
first set of stable identifications are shown (dashed-dotted
lines). Further stable identifications are obtained from these by
iterating the map backwards. The points of the period 2 orbit of
infinitesimal triangles give rise to a period 2 orbit of 3-pronged
singularities. There are infinitely many 1-pronged singularities (two
of which come from the points marked with $\bullet$ in the figure)
converging in the future to the period 6 orbit (indicated by $\circ$)
corresponding to the period 6 attracting orbit of $F_2$ and converging
in the past to the period 3 orbit at infinity (indicated by
\rule{1.5mm}{1.5mm}). There is also an orbit of 3-pronged
singularities with the same forward and backward fates as the orbit of
1-pronged singularities. A detailed view of the shaded regions in
Figure~\ref{fig:thtrident} is shown in Figure~\ref{fig:bubbleident}.
\begin{figure}
\pichere{1}{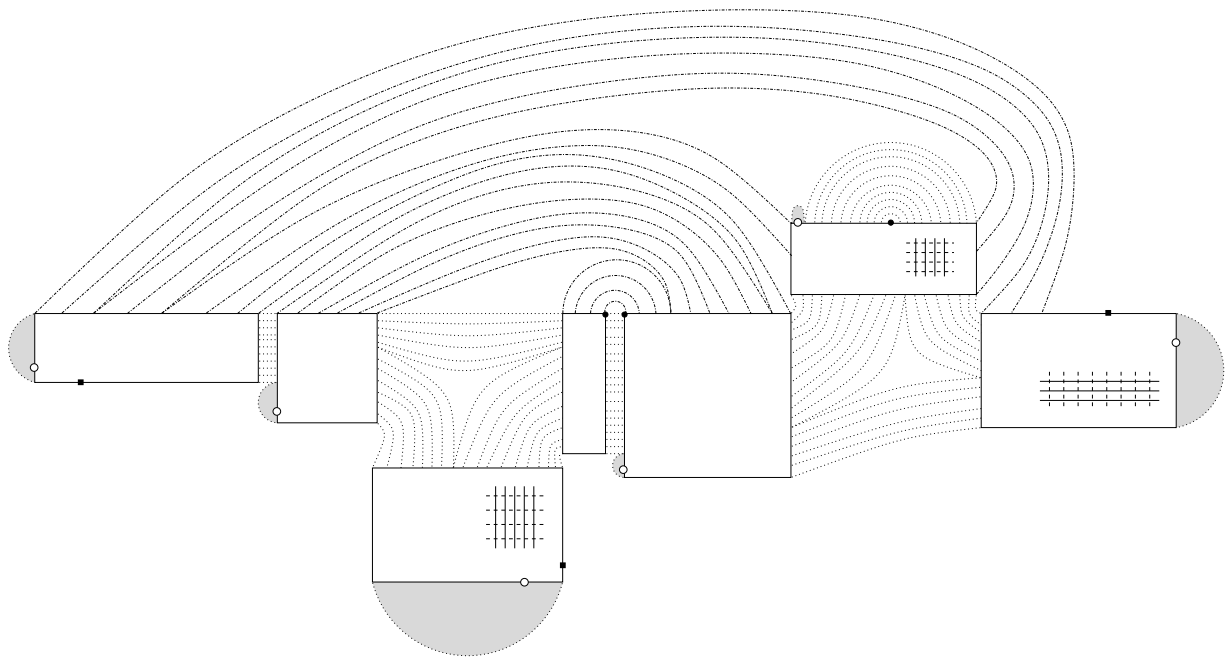}
\caption{The identifications for $F_2$.}
\label{fig:thtrident}
\end{figure}
\begin{figure}
\pichere{0.4}{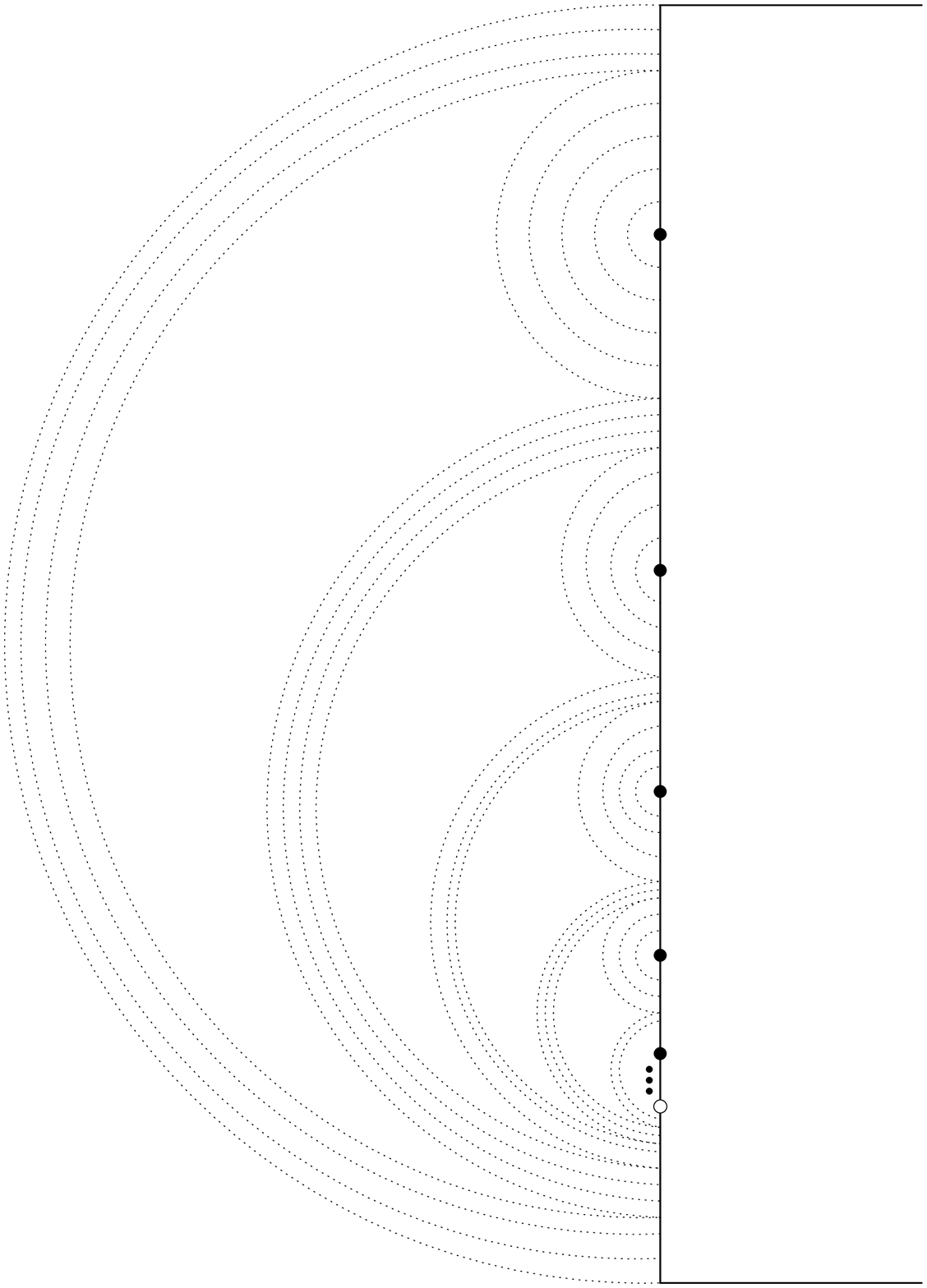}
\caption{The identifications in the shaded regions in
Figure~\ref{fig:thtrident}.} 
\label{fig:bubbleident}
\end{figure}
\end{example}

\subsection{1-pronged singularities}
\label{subsec:1prongs}

The generalized pseudo-Anosov maps constructed in the previous section
preserve a pair of measured foliations with singularities. The
1-pronged singularities are the most important in many interrelated
ways.  Dynamically, they play a role analogous to that played by the
critical points for endomorphisms of the interval. They also `hold the
map in place,' so to say, in the sense that isotopies relative to the
set of 1-pronged singularities cannot destroy any
dynamics.\footnote{This has been proved for the case of finitely many
singularities~\cite{Hall,Handel,Thurston} and is conjectured to be
true in general~\cite{Hall-Boyland}.} 

We now describe how to find the orbits of 1-pronged singularities of
the invariant foliations $\cF^{u,s}$. They come from infinitesimal
1-gons of $\tau$ so we need to be able to determine the orbits of
these. It follows from Proposition~\ref{prop:n-gons} that if a 1-gon
contains a periodic point, it is itself periodic under $\phi$ and
therefore corresponds to a periodic 1-pronged singularity of the
invariant foliations. If a 1-gon does not belong to a periodic orbit,
then there exists a {\em first} 1-gon whose orbit contains the given
1-gon, that is, there exists a 1-gon $\Delta$, which is not the image
of any 1-gon under $\phi$ and whose orbit contains the given 1-gon.
In this case, there exists either a real edge $e$ such that
$\phi(e)\supset\Delta$ or there is an infinitesimal edge $e'$, which
is not a 1-gon, that is, whose boundary points are distinct switches,
such that $\phi(e')=\Delta$. In either case, there must exist a real
edge $e$, an arc $\gamma\subset e$ and a smallest integer $k\ge 1$
such that $\phi^k(\gamma)=\Delta$. If there are several such arcs
(possibly contained in several distinct real edges), we use the
ambient surface or thick graph to choose $\gamma$ to be {\em
innermost} among of them. By this we mean the following. To the arc
$\gamma$ there corresponds a {\em thick arc}
$\Gamma=[a,b]\times[0,1]\subset s$ where $s$ is a strip of
$\ofg$. Because $\phi^k(\gamma)=\Delta$, $F^n(\Gamma)$ is contained in
the junction $J$ that contains $\Delta$ and $F^n(\{a,b\}\times[0,1])$
are contained in the same component of $\partial J$, that is,
$F^n(\Gamma)$ `makes a turn' inside $J$. We call $\gamma$ innermost if
$F^n(\Gamma)$ is an innermost turn among all thick arcs that map to
$J$ under $F^n$. Since infinitesimal edges are assigned 0 length in
the construction of the invariant foliations, $\gamma$ corresponds to
a vertical segment in the rectangle $R_e$ associated to the real edge
$e$. Under the identifications required to make the leaves of the
stable foliation infinite, $R_e$ will be folded and one of the
endpoints of this segment will become a 1-pronged singularity (the one
that becomes innermost after folding). Since the `horizontal' sides
of rectangles are contained in the unstable manifold of the point at
infinity, it follows that this 1-pronged singularity is backward
asymptotic to the point at infinity. To summarize, we have proved the
\begin{thm}
Let $\cF^{u,s}$ be the stable and unstable foliations of the
generalized pseudo-Anosov map associated to a Markov thick graph map
as constructed above.  Then the 1-pronged singularities of $\cF^{u,s}$
either belong to periodic orbits or are backward asymptotic to the
point at infinity and forward asymptotic to one of finitely many
periodic orbits.
\end{thm}

There is an easy way of tracking the orbits of 1-pronged singularities
backwards using the graph map \gmx. Suppose $p\in G$ is a critical
point of $f$ (that is, a point at which $f$ is not a local
homeomorphism), so that $v=f(p)$ is a vertex. Assume, without loss of
generality, that $f(p)$ is innermost in the junction
$V=\pi^{-1}(v)$. Choose $p_{-1}\in f^{-1}(p)$ so that $f^2(p_{-1})$ is
innermost in $V$, then choose $p_{-2}\in f\I(p_{-1})$ so that
$f^3(p_{-2})$ is innermost in $V$ and proceed like this. If we reach a
point $p_{-n}$ which is not a vertex, then all subsequent ones are not
vertices (since vertices map to vertices) and, in fact, at each
subsequent step, there is exactly one innermost preimage (since the
graph map is assumed to be onto). In this case, $p_{-n}$ converges to
the orbit at infinity. Otherwise, there are two possibilities: either
we eventually return to $p$ or we get to a point which has {\em no}
innermost preimage. In the first case, we have found a periodic orbit
of 1-pronged singularities and in the second, we were not following an
orbit of 1-pronged singularities. Below we give examples explaining
all three possibilities.

\begin{examples}
\begin{enumerate}[a)]
\item Consider the thick tree map $F_2$. In Figure~\ref{fig:trmap2},
the quotient tree map is drawn. To avoid having too many symbols in
the figure, the images of points in $T$ under $f$ are denoted by the
same symbol on the image tree $f(T)$ drawn below. The only
critical point is $p=4$ and $v=f(4)=5$. The preimage $f\I(4)$ has two
points: 2 and a point in edge $e_5$. It is the latter that becomes
innermost under $f^2$, so $p_{-1} \in e_5$. The preimage $f\I(p_{-1})$
again consists of two points, of which the one that becomes innermost
under $f^3$ lies on edge $e_2$: this is $p_{-2}$. Again $f\I(p_{-2})$
consists of two points and $p_{-3}\in e_1$. From here on things
repeat, as shown in Figure~\ref{fig:trmap2}: $p_{-4}\in e_7$,
$p_{-5}\in e_3$, $p_{-6}\in e_1$, $p_{-7}\in e_7$, etc. The sequence
$p_n$ converges to the period 3 orbit at infinity.
\begin{figure}
\lab{f}{f_2}{}
\lab{a}{e_1}{}
\lab{b}{e_2}{}
\lab{c}{e_3}{}\lab{d}{e_4}{}\lab{e}{e_5}{}\lab{g}{e_6}{}
\lab{i}{e_7}{}
\lab{j}{f(e_1)}{r}\lab{k}{f(e_2)}{}
\lab{l}{f(e_3)}{}
\lab{m}{f(e_4)}{bl}
\lab{n}{f(e_5)}{}\lab{p}{f(e_7)}{b}
\lab{1}{p_{-1}}{}\lab{2}{p_{-2}}{}\lab{3}{p_{-3}}{}\lab{4}{p_{-4}}{b}
\lab{5}{p_{-5}}{l}\lab{6}{p_{-6}}{l}\lab{7}{p_{-7}}{b}\lab{8}{p_{-8}}{l}
\lab{9}{p_{-9}}{l}
\lab{x}{4}{}
\pichere{0.6}{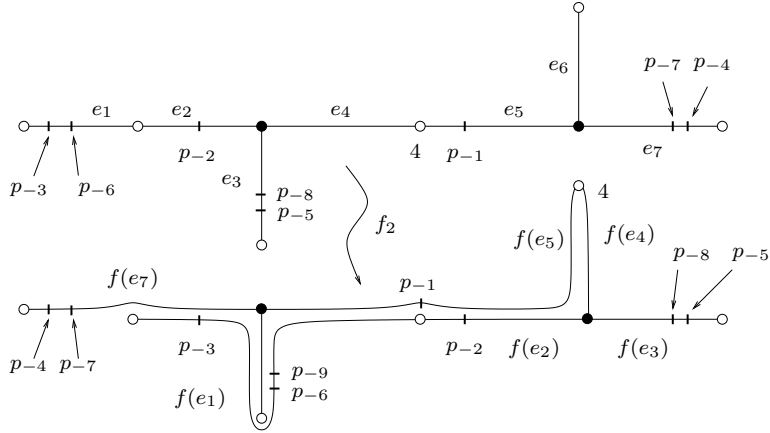}
\caption{The backward orbit of 1-pronged singularities for $F_2$.}
\label{fig:trmap2}
\end{figure}    
\item In this example (see Figure~\ref{fig:10010}, where only the
induced tree endomorphism is drawn) $p=3$ and $v=5$. The preimage
$f\I(p)$ has two points, namely, $4$ and a point on edge $e_1$. Under
$f^2$ it is point 4 that becomes innermost so 
$p_{-1}=4$. Again, $f\I(4)$ contains two points, namely, 2 and a point
on edge $e_3$, but it is 2 that becomes innermost under $f^3$ so that
$p_{-2}=2$. Continuing like this, we get $p_{-3}=1$, $p_{-4}=5$ and
$p_{-5}=3$. Thus, $p$ lies gives rise to a periodic orbit of 1-pronged
singularities. The associated generalized pseudo-Anosov is, in fact, a
pseudo-Anosov map in the sense of Thurston.
\begin{figure}
\lab{1}{1}{bl}\lab{2}{2}{bl}\lab{3}{3}{bl}\lab{4}{4}{bl}\lab{5}{5}{bl}
\lab{a}{e_1}{}
\lab{b}{e_2}{}
\lab{c}{e_3}{}\lab{d}{e_4}{}\lab{e}{f(e_1)}{}\lab{g}{f(e_2)}{}
\lab{h}{f(e_3)}{b}
\lab{i}{f(e_4)}{b}\lab{f}{f}{l}
\pichere{0.4}{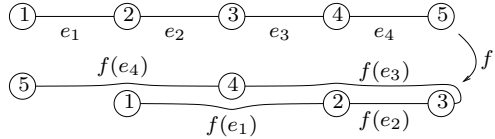}
\caption{The quotient of a thick tree map whose associated
generalized pseudo-Anosov has a periodic orbit of 1-pronged
singularities.}
\label{fig:10010}
\end{figure}    

\item The third possibility is shown in Figure~\ref{fig:no1prong}. The
top part of the figure contains parts of the graph containing the
points $p$ and $v=f(p)$ and two points in $f\I(p)$; the bottom
contains their images under $f$. Because no arc has image entering $p$
through edge $a$ and exiting through edge $b$ (indicated by the dashed
line in the figure), $f\I(p)$ does not
contain a point whose image under $f^2$ is innermost in $v$.
\begin{figure}
\lab{f}{f}{}\lab{p}{p}{}\lab{v}{v}{}\lab{x}{f\I(p)}{}\lab{y}{f(p)=v}{l}
\lab{a}{a}{b}\lab{b}{b}{}\lab{c}{f(a)}{b}\lab{d}{f(b)}{}
\pichere{0.4}{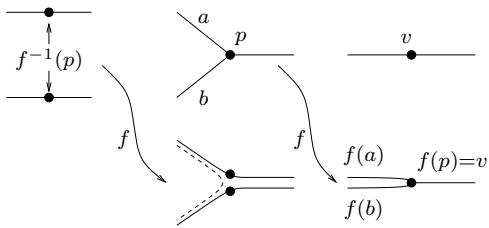}
\caption{When a critical point of a tree map does not produce an orbit
of 1-pronged singularities.}
\label{fig:no1prong} 
\end{figure}
\end{enumerate}
\end{examples} 

\subsection{The complex structure}

We now describe a natural complex structure associated to a
generalized pseudo-Anosov map, with respect to which it becomes a
Teichm\"uller mapping and the invariant foliations become the
horizontal and vertical trajectories of an associated integrable quadratic
differential.

As was seen above, $S$ is obtained as the quotient of Euclidean
rectangles under side identifications by Euclidean isometries. In the
interior of the rectangles, the complex structure is the one
determined by the Euclidean structure. At $k$-pronged singularities,
use the maps $z\mapsto z^{2/k}$ to define coordinate charts. This
produces a complex structure at all points of $S$, except at the
accumulations of singularities. These are topologically isolated and
it is necessary to decide whether the complex structure in the complement
regards them as punctures or as holes. We now argue that they are in
fact punctures so that the complex structure extends across them
uniquely. To prove this we present a sequence of concentric annuli
converging to an accumulation of singularities and whose moduli add up
to infinity. By length-area arguments (see~\cite{Ah,LV}), the result
follows\footnote{I am grateful to Fred Gardiner for suggesting this argument
(see also~\cite{GaEa}).}. 

The basic inequality used is 
\[ \mod(A)\ge \frac{\inf_{\gamma\in\Gamma}{l(\gamma)}^2}{\mbox{Area}(A)} \]
where $\Gamma$ is the set of all rectifyable curves joining the
boundary components of the annular region $A$ and $l(\gamma)$ is the
length of the curve $\gamma$. 

We first give the argument for the accumulation of singularities $p$ in
the generalized pseudo-Anosov associated to the horseshoe.

\begin{figure}
\pichere{0.4}{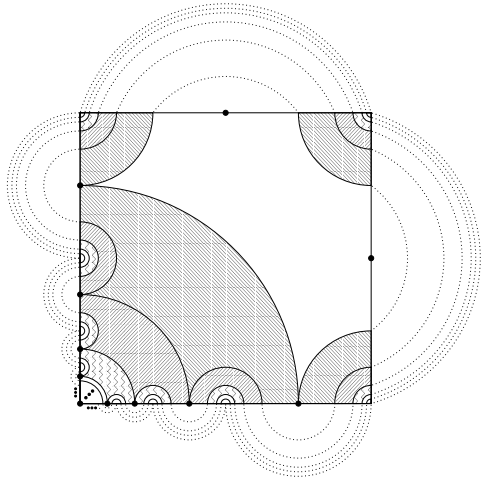}
\caption{A sequence of concentric annular regions whose sum of moduli
diverges.} 
\label{fig:infty}
\end{figure}

In Figure~\ref{fig:infty} are shown the first four concentric annular
regions of a sequence $\{A_n\}_{n=1}^\infty$ converging to $p$. We
argue that $\mod(A_n)\propto 1/n$ (that is, $\mod(A_n)$ is
proportional to $1/n$), implying thus that
$\sum_{n=1}^\infty\mod(A_n)$ diverges. Annular region $A_n$ is made up
of three round quarter-annuli at the upper-left and right and
lower-right corners, the region bounded by four circular arcs on the
lower-left and $2n$ round half-annuli, $n$ along each of the left and
bottom sides of the square. All parts being (parts of) Euclidean
annuli, they contribute to the extremal length of the family of all
curves joining the boundary components of $A_n$ at least like
Euclidean annuli do. The infimum of the lengths of curves joining
boundary components is $\propto 1/2^n$. The sum of the areas of the
quarter annuli is $\propto 1/2^{2n}$, whereas the $2n$ half-annuli
contribute $\propto n/2^{2n}$ to the total area. Therefore, we have
\begin{eqnarray*}
\mod(A_n) & \ge & \frac{(\mbox{Distance between boundary components
of\ }A_n)^2}{\mbox{Area}(A_n)} \\
          & \ge & \frac{C_1}{C_2 + C_3 n} 
\end{eqnarray*}

The general result is 

\begin{thm}
Let $F\colon(S,\ofg)\to(S,\ofg)$ be a Markov thick graph map whose associated
transition matrix is irreducible and aperiodic. Let $\Phi\colon S\to
S$ be the associated generalized pseudo-Anosov and
$\{p_1,\ldots,p_k\}$ the accumulation points of pronged singularities
of the $\Phi$-invariant foliations.  Then the complex structure on
$S\setminus\{p_1,\ldots,p_k\}$ induced by the Euclidean structure
extends uniquely to a complex structure on the compact surface $S$.
\end{thm}

\medskip
\noindent
{\em Sketch of Proof.}
Let $p$ be an accumulation of pronged singularities of the invariant
foliations corresponding to an attracting periodic orbit of $F$, which
we will also denote by $p$ (the other possibility, namely, that $p$ is
the point at infinity, is treated similarly). Passing to an appropriate
iterate, we may assume that $p$ is an attracting fixed point of
$F$. For simplicity, we assume it lies in a 1-junction $V$; the other
cases are analogous.

\begin{figure}
\pichere{0.4}{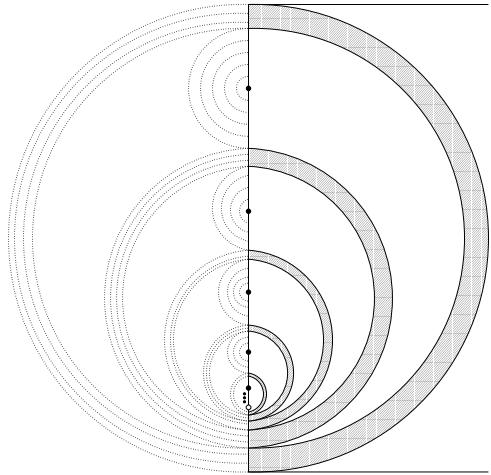}
\caption{Another sequence of concentric annular regions whose sum of moduli
diverges.} 
\label{fig:puncture}
\end{figure}

There are two possibilities to consider: either the invariant train
track $\tau$ has a bubble in $V$ which encloses $p$ or not. If there
is a bubble $\beta$ enclosing $p$, all its iterates also do. They then
give rise to an infinite collection of disjoint concentric annular regions
$A_n$, all enclosing $p$, with $\mod(A_n)\propto 1$ (see
Figure~\ref{fig:puncture}). In case there is no bubble enclosing $p$
in $V$, then $p$ is the point at infinity and the situation is
analogous to the horseshoe example explained above.  \qed

\bibliography{infinitett}

\end{document}